\documentclass[11pt]{article}
\usepackage{graphicx}
\usepackage{amsfonts}
\usepackage{theorem}
\theorembodyfont{\rmfamily}
\newtheorem{Lem}{Lemma}[section]

\newtheorem{The}[Lem]{Theorem}
\newtheorem{Prop}[Lem]{Proposition}

\newtheorem{Ex}[Lem]{Example}

\newtheorem{Rem}[Lem]{Remark}

\newtheorem{Con}[Lem]{Conjecture}
\newcommand{\qed}{\hbox{\rule{6pt}{6pt}}}

\begin{document}
\title{Trace inequalities for products of matrices}
\author{Shigeru Furuichi$^1$\footnote{E-mail:furuichi@cssa.chs.nihon-u.ac.jp},  Ken Kuriyama$^2$\footnote{E-mail:kuriyama@yamaguchi-u.ac.jp} and Kenjiro Yanagi$^2$\footnote{E-mail:yanagi@yamaguchi-u.ac.jp}\\
$^1${\small Department of Computer Science and System Analysis,}\\
{\small College of Humanities and Sciences, Nihon University,}\\
{\small 3-25-40, Sakurajyousui, Setagaya-ku, Tokyo, 156-8550, Japan}\\
$^2${\small Division of Applied Mathematical Science,}\\
{\small Graduate School of Science and Engineering, Yamaguchi University,}\\
{\small 2-16-1, Tokiwadai, Ube City, 755-8611, Japan}}
\date{}
\maketitle
{\bf Abstract.} In this short paper, we study some trace inequalities of the products of the matrices and the power of matrices by the use of 
elementary calculations.

\vspace{3mm}

\vspace{3mm}
\vspace{3mm}

\section{Introduction}
In this short paper, we treat some matrix trace inequalities.
Let $M(n,\mathbb{C})$ be the set of all $n\times n$ matrices on the complex field $\mathbb{C}$.
For a Hermitian matrix $X$, $X \geq 0$ means that $\langle \phi\vert  X \vert \phi \rangle \geq 0$ for any vector $\vert \phi \rangle \in \mathbb{C}^n$.
Here we denote the set of all Hermitian matrices by $M_h(n,\mathbb{C})$, that is, $M_h(n,\mathbb{C})\equiv \left\{X\in M(n,\mathbb{C}) \vert X=X^* \right\}$.
In addition, we denote the set of all nonnegative matrices by $M_+(n,\mathbb{C}) $, that is,  $M_+(n,\mathbb{C}) \equiv \left\{X\in M(n,\mathbb{C}) \vert X\geq 0 \right\}$.
Then we have the following propositions.

\begin{Prop}
For two real valued functions $f,g$ on $D \subset \mathbb{R}$ and $L \in M_+(n,\mathbb{C}), A \in M_h(n,\mathbb{C})$, we have,
\begin{equation}\label{intro_ineq01}
Tr\left[f(A)Lg(A)L\right] \leq \frac{1}{2} Tr\left[(f(A)L)^2 + (g(A)L)^2\right].  
\end{equation}
\end{Prop}
{\it Proof:}
Since the arithmetic mean is greater than the geometric mean and by the use of Schwarz's inequality,
we have
\begin{eqnarray*}
Tr\left[f(A)Lg(A)L\right] &=& Tr\left[(L^{1/2}f(A)L^{1/2})(L^{1/2}g(A)L^{1/2})\right] \\
&\leq& \left(Tr\left[(f(A)L)^2\right]\right)^{1/2}  \left(Tr\left[(g(A)L)^2\right]\right)^{1/2} \\ 
&\leq& \frac{1}{2} Tr\left[(f(A)L)^2 + (g(A)L)^2\right].  
\end{eqnarray*}
\hfill \qed

\begin{Prop}
For two real valued functions $f,g$ on $D \subset \mathbb{R}$ and  $L, A \in M_h(n,\mathbb{C})$, we have,
\begin{equation}\label{intro_ineq02}
Tr\left[f(A)Lg(A)L\right] \leq \frac{1}{2} Tr\left[f(A)^2L^2 + g(A)^2L^2\right]. 
\end{equation}
\end{Prop}
{\it Proof:}
Since the arithmetic mean is greater than the geometric mean and by the use of Schwarz's inequality,
we have
\begin{eqnarray*}
Tr\left[f(A)Lg(A)L\right] &\leq& \left(Tr\left[Lf(A)^2L\right]\right)^{1/2}  \left(Tr\left[Lg(A)^2L\right]\right)^{1/2} \\
&\leq& \frac{1}{2} Tr\left[f(A)^2L^2 + g(A)^2L^2\right].
\end{eqnarray*}
\hfill \qed

Note that the right hand side of the inequalities (\ref{intro_ineq01}) is less than the  right hand side of the inequalities (\ref{intro_ineq02}),
since we have $Tr\left[XYXY\right] \leq Tr\left[X^2Y^2\right]$ for two Hermitian matrices $X$ and $Y$ in general. 

\begin{The}
For two real valued functions $f,g$ on $D \subset \mathbb{R}$ and $L, A \in M_h(n,\mathbb{C})$, if we have
$f(x) \leq g(x)$ or $f(x) \geq g(x)$ for any $x \in D$,
then we have,
\begin{equation}
Tr\left[f(A)Lg(A)L\right]  \leq \frac{1}{2}  Tr\left[(f(A)L)^2 + (g(A)L)^2\right].
\end{equation}
\end{The}
{\it Proof:}
For a spectral decomposition of $A$ such as $A = \sum_k \lambda_k \vert \phi_k \rangle \langle \phi_k \vert$, we have
$$
Tr\left[f(A)Lg(A)L\right] = \frac{1}{2} \sum_{m,n} \left\{  f(\lambda_m)g(\lambda_n)+g(\lambda_m)f(\lambda_n) \right\} \vert \langle \phi_m\vert L\vert \phi_n \rangle \vert^2
$$
and
$$
 \frac{1}{2}  Tr\left[(f(A)L)^2 + (g(A)L)^2\right] = \frac{1}{2} \sum_{m,n} \left\{  f(\lambda_m)f(\lambda_n)+g(\lambda_m)g(\lambda_n) \right\} \vert \langle \phi_m\vert L\vert \phi_n \rangle \vert^2.
$$
Thus we have the present theorem by
$$
f(\lambda_m)f(\lambda_n)+g(\lambda_m)g(\lambda_n) - f(\lambda_m)g(\lambda_n)-g(\lambda_m)f(\lambda_n) = 
\left\{f(\lambda_m) - g(\lambda_m) \right\}\left\{f(\lambda_n) - g(\lambda_n) \right\} \geq 0.
$$
\hfill \qed

Trace inequalities for multiple products of two matrices have been studied by Ando,Hiai and Okubo in \cite{AHO} with the notion of majorization.
Our results in the present paper are derived by the elementary calculations without the notion of majorization.

Next we consider the further specialized forms for the products of matrices.
We have the following trace inequalities on the products of the power of matrices.
\begin{Prop} \label{prop01}
\begin{itemize}
\item[(i)] For any natural number $m$ and $T,A \in M_+(n,\mathbb{C})$, we have the inequality
\begin{equation} \label{ineq01}
Tr\left[ \left(T^{1/m}A\right)^m\right] \leq Tr\left[ TA^m\right].
\end{equation}
\item[(ii)] For $\alpha \in [0,1]$, $T \in M_+(n,\mathbb{C})$ and $A \in M_h(n,\mathbb{C})$, we have the inequalities
\begin{equation} \label{ineq02}
Tr\left[ \left(T^{1/2}A\right)^2\right] \leq Tr\left[T^{\alpha}AT^{1-\alpha}A\right] \leq Tr\left[ TA^2\right].
\end{equation}
\end{itemize}
\end{Prop}
{\it Proof}:
\begin{itemize}
\item[(i)] Putting $p=m,r=1/m,X=A^m$ and $Y=T$ in Araki's inequality \cite{Araki}:
$$
Tr \left[ \left(Y^{r/2}X^rY^{r/2}\right)^p  \right]\leq Tr \left[\left(Y^{1/2}XY^{1/2}\right)^{rp}  \right]
$$
for $X,Y \in M_+(n,\mathbb{C})$ and $p>0,0\leq r \leq 1$, we have the inequality (\ref{ineq01}).
\item[(ii)] By the use of Lemma \ref{lem01} in the below, we straightforwardly have 
$$Tr\left[T^{\alpha}AT^{1-\alpha}A\right] \leq Tr\left[ TA^2\right].$$
Again by the use of Lemma \ref{lem01}, we have
\begin{eqnarray*}
Tr\left[ {T^{1/2} AT^{1/2} A} \right] &=& Tr\left[ {\left( {T^{1/4} AT^{1/4} } \right)^2 } \right] \\
 &\le& Tr\left[ {T^{\alpha  - 1/2} \left( {T^{1/4} AT^{1/4} } \right)T^{1/2 - \alpha } \left( {T^{1/4} AT^{1/4} } \right)} \right] \\
 &=& Tr\left[ {T^\alpha  AT^{1 - \alpha } A} \right].
\end{eqnarray*}
\end{itemize}
\hfill \qed

\begin{Lem}  {\bf (Bourin \cite{Bourin}, Fujii \cite{Fujii})} \label{lem01}
Let $f$ and $g$ be functions on the domain $D \subset \mathbb{R}$.
For any Hermitian matrices $A$ and $X$ having spectrums in $D$, we have
\begin{itemize}
\item[(i)] If $(f,g)$ satisfies $(f(a) - f(b)) (g(a) - g(b)) \geq 0$ for any $a,b \in D$, then 
$$Tr[f(A)Xg(A)X] \leq Tr[f(A)g(A)X^2]. $$
\item[(ii)] If $(f,g)$ satisfies $(f(a) - f(b)) (g(a) - g(b)) \leq 0$ for any $a,b \in D$, then 
$$Tr[f(A)Xg(A)X] \geq Tr[f(A)g(A)X^2]. $$
\end{itemize}
\end{Lem}

\section{Main results}
In this short note, we study a generalization of Proposition \ref{prop01}.
To this end, we prepare the elementary inequality on the arithmetic mean and the geometric mean. See \cite{HLP} for example.
\begin{Lem} \label{lem02} 
For positive numbers $a_1,a_2,\cdots,a_m$ and $p_1,p_2,\cdots,p_m$ with $p_1+p_2+\cdots +p_m=1$, the following inequality holds:
$$ a_1^{p_1}a_2^{p_2}\cdot \dots \cdot a_m^{p_m} \leq p_1 a_1 +p_2 a_2 +\cdots +p_m a_m,$$
with equality if and only if $a_1=a_2=\cdots =a_m$.
\end{Lem}

\begin{The}   \label{the01}
For positive numbers $p_1,p_2,\cdots,p_m$ with $p_1+p_2+\cdots +p_m=1$ and $T,A\in M_+(2,\mathbb{C})$, we have the inequalities
\begin{equation} \label{ineq03}
Tr\left[ \left( T^{1/m}A\right)^m \right]\leq Tr\left[ T^{p_1}AT^{p_2}A\cdots T^{p_m}A \right] \leq Tr\left[ TA^m \right].
\end{equation}
\end{The}
{\it Proof}:
We write $T  = \sum\limits_{i = 0}^1 {\lambda _i \left| {\psi _i } \right\rangle \left\langle {\psi _i } \right|} $, where 
we can take $\left\{{\left| {\psi _0 } \right\rangle },{\left| {\psi _1 } \right\rangle }\right\}$ as a complete orthonormal base.

By the use of Lemma \ref{lem02}, we then have
\begin{eqnarray}
 && Tr\left[ {T ^{p_1 } AT ^{p _2 } A \cdots T ^{p _m } A} \right] \nonumber \\ 
 && = \sum\limits_{i_1 ,i_2 , \cdots ,i_m } {\left( {\frac{{\lambda _{i_1 }^{p _1 } \lambda _{i_2 }^{p _2 }  \cdots \lambda _{i_m }^{p _m }  + \lambda _{i_1 }^{p _2 } \lambda _{i_2 }^{p _3 }  \cdots \lambda _{i_m }^{p _1 }  +  \cdots  + \lambda _{i_1 }^{p _m } \lambda _{i_2 }^{p _1 }  \cdots \lambda _{i_m }^{p _{m - 1} } }}{m}} \right)} \nonumber \\
 && \hspace*{1cm} \times \left\langle {\psi _{i_1 } } \right|A\left| {\psi _{i_2 } } \right\rangle \left\langle {\psi _{i_2 } } \right|A\left| {\psi _{i_3 } } \right\rangle  \cdots \left\langle {\psi _{i_m } } \right|A\left| {\psi _{i_1 } } \right\rangle \nonumber  \\ 
 && \le \sum\limits_{i_1 ,i_2 , \cdots ,i_m } {\left( {\frac{{p _1 \lambda _{i_1 }   +  \cdots  + p _m \lambda _{i_m }  + p_2\lambda_{i_1}+\cdots +p_1\lambda_{i_m}+\cdots + p _m \lambda _{i_1 }  +  \cdots  + p _{m - 1} \lambda _{i_m } }}{m}} \right)}\nonumber  \\
 && \hspace*{1cm} \times \left\langle {\psi _{i_1 } } \right|A\left| {\psi _{i_2 } } \right\rangle \left\langle {\psi _{i_2 } } \right|A\left| {\psi _{i_3 } } \right\rangle  \cdots \left\langle {\psi _{i_m } } \right|A\left| {\psi _{i_1 } } \right\rangle  \nonumber \\ 
 && = \sum\limits_{i_1 ,i_2 , \cdots ,i_m } {\left( {\frac{{\lambda _{i_1 }  + \lambda _{i_2 }  +  \cdots  + \lambda _{i_m } }}{m}} \right)} \left\langle {\psi _{i_1 } } \right|A\left| {\psi _{i_2 } } \right\rangle \left\langle {\psi _{i_2 } } \right|A\left| {\psi _{i_3 } } \right\rangle  \cdots \left\langle {\psi _{i_m } } \right|A\left| {\psi _{i_1 } } \right\rangle \nonumber  \\ 
 && = Tr\left[ {T A^m } \right]. \label{the_ineq01}
 \end{eqnarray}
We should note that the above calculation is assured by
$$\left\langle {\psi _{i_1 } } \right|A\left| {\psi _{i_2 } } \right\rangle \left\langle {\psi _{i_2 } } \right|A\left| {\psi _{i_3 } } \right\rangle  \cdots \left\langle {\psi _{i_m } } \right|A\left| {\psi _{i_1 } } \right\rangle \geq 0,$$
since every $i_j\,(j=1,2,\cdots ,m)$ takes $0$ or $1$. See Lemma \ref{lem03} in the below.

Again by the use of Lemma \ref{lem02} with $p_i=1/m$, we have
\begin{eqnarray}
&& Tr\left[ {T ^{p_1 } AT ^{p_2 } A \cdots T ^{p_m } A} \right] \nonumber \\ 
&&  = \sum\limits_{i_1 ,i_2 , \cdots ,i_m } {\left( {\frac{{\lambda _{i_1 }^{p_1 } \lambda _{i_2 }^{p_2 }  \cdots \lambda _{i_m }^{p_m }  + \lambda _{i_1 }^{p_2 } \lambda _{i_2 }^{p_3 }  \cdots \lambda _{i_m }^{p_1 }  +  \cdots  + \lambda _{i_1 }^{p_m } \lambda _{i_2 }^{p_1 }  \cdots \lambda _{i_m }^{p_{m - 1} } }}{m}} \right)}\nonumber \\
&& \hspace*{1cm} \times   \left\langle {\psi _{i_1 } } \right|A\left| {\psi _{i_2 } } \right\rangle \left\langle {\psi _{i_2 } } \right|A\left| {\psi _{i_3 } } \right\rangle  \cdots \left\langle {\psi _{i_m } } \right|A\left| {\psi _{i_1 } } \right\rangle  \nonumber \\ 
&& \ge \sum\limits_{i_1 ,i_2 , \cdots ,i_m } {\lambda _{i_1 }^{1/m} \lambda _{i_2 }^{1/m}  \cdots \lambda _{i_m }^{1/m} \left\langle {\psi _{i_1 } } \right|A\left| {\psi _{i_2 } } \right\rangle \left\langle {\psi _{i_2 } } \right|A\left| {\psi _{i_3 } } \right\rangle  \cdots } \left\langle {\psi _{i_m } } \right|A\left| {\psi _{i_1 } } \right\rangle  \nonumber \\ 
&&  = Tr\left[ {\left( {T ^{1/m} A} \right)^m } \right]. \label{the_ineq02}
\end{eqnarray}
Note that the inequalities (\ref{the_ineq01}) and (\ref{the_ineq02}) hold even if $\lambda_0=0$ or $\lambda_1=0$. It is a trivial when $\lambda_0=0$ and $\lambda_1=0$. 
Thus the proof of the present theorem is completed.

\hfill \qed

Note that the second inequality of (\ref{ineq03}) is derived by putting $f_i(x)=x^{p_i}$ and $g_i(x)=x$ for $i=1,\cdots,m$ in Theorem 4.1 of \cite{AHO}.
However the first inequality of (\ref{ineq03}) can not be derived by applying Theorem 4.1 of \cite{AHO}.

\begin{Rem}
From the process of the proof in Theorem \ref{the01}, we find that,
if $T$ is an invertible, then the equalities in both inequalities (\ref{the_ineq01}) and (\ref{the_ineq02}) hold if and only if $T=kI$.
\end{Rem}

\begin{Lem} \label{lem03}
For $A\in M_+(2,\mathbb{C})$ and a complete orthonormal base $\left\{\psi_0,\psi_1\right\}$ of $\mathbb{C}^2$, we have
 $$\left\langle {\psi _{i_1 } } \right|A\left| {\psi _{i_2 } } \right\rangle \left\langle {\psi _{i_2 } } \right|A\left| {\psi _{i_3 } } \right\rangle  \cdots \left\langle {\psi _{i_m } } \right|A\left| {\psi _{i_1 } } \right\rangle \geq 0,$$
where every $i_j\,(j=1,2,\cdots ,m)$ takes $0$ or $1$. 
\end{Lem}
{\it Proof}:
We set a symmetric group by
\[
\pi  \equiv \,\,\left( \begin{array}{l}
 \,\,1\,\,\,\,\,\,\,2\,\,\,\, \cdots \,\,\,\,\,\,\,m \\ 
 \,m\,\,\,\,\,\,1\,\,\,\, \cdots \,\,\,\,m - 1 \\ 
 \end{array} \right).
\]
We also set 
\[
S \equiv \left\{ {1 \le j \le m\left| {i_j  = i_{\pi \left( j \right)} } \right.} \right\}
\]
for $\left(i_1,i_2,\cdots,i_m\right) \in \left\{0,1\right\}^m$.
Then we have
\[
\prod\limits_{j = 1}^m {\left\langle {\psi _{i_{\pi  \left( j \right)} } } \right|A\left| {\psi _{i_j } } \right\rangle }  = \prod\limits_{j \in S} {\left\langle {\psi _{i_{\pi  \left( j \right)} } } \right|A\left| {\psi _{i_j } } \right\rangle }  \cdot \prod\limits_{j \notin S} {\left\langle {\psi _{i_{\pi  \left( j \right)} } } \right|A\left| {\psi _{i_j } } \right\rangle }. 
\]
Here we have
\[
\prod\limits_{j \in S} {\left\langle {\psi _{i_{\pi  \left( j \right)} } } \right|A\left| {\psi _{i_j } } \right\rangle }  \ge 0,
\]
since $A$ is a nonnegative matrix.
In addition, $m-\vert S\vert$ necessarily takes $0$ or an even number (see Lemma \ref{lem04} in the below) and then we have
\[
\prod\limits_{j \notin S} {\left\langle {\psi _{i_{\pi  \left( j \right)} } } \right|A\left| {\psi _{i_j } } \right\rangle }  
=\left| {\left\langle {\psi _0 } \right|A\left| {\psi _1 } \right\rangle } \right|^{m - \left| S \right|} \ge 0.
\]
Therefore we have the present lemma.

\hfill \qed

\begin{Lem} \label{lem04}
Image the situation that we put arbitrary $l_0$ vectors $\vert \psi_0 \rangle $ and
 $l_1$ vectors $\vert \psi_1 \rangle$ on the circle 
and then the circle is divided into $l_0+l_1$ circular arcs.
Then the number $\vert S^c\vert$ of the circular arcs having different edges is $0$ 
or an even number.
\end{Lem}
{\it Proof}:
If $l_0=0$ or $l_1=0$, then  $\vert S^c\vert=0$. Thus we consider $l_0\neq 0$ and $l_1\neq 0$.
We suppose that the number of the circular arcs having same $\vert \psi_0 \rangle $ in both edges is $E$,
and the number of the circular arcs having $\vert \psi_0 \rangle $ and $\vert \psi_1 \rangle $ in their edges is $F$.
We now form the circular arcs by combining every $\vert \psi_0 \rangle $ with its both sides. (We do not consider the circular arcs formed by the other method.)
Thus the number of the circular arcs formed by the above method is an even number, since every $\vert \psi_0 \rangle $ forms two circular arcs. 
In addition, its number coincides
with $2E+F$ which shows the number that every $\vert \psi_0 \rangle $ is doubly counted. Thus $2E+F$ takes an even number.
Then $F$ must be an even number.

\hfill \qed

We give an example of Lemma \ref{lem03} for readers' convenience.

\begin{Ex}
For $\left\{i_1,i_2,i_3,i_4,i_5,i_6,i_7 \right\}=\left\{0,0,1,1,0,1,0\right\}$,
we have $S=\left\{1,2,4\right\}$ and $S^c=\left\{3,5,6,7\right\}$.
Then we have
\[
\prod\limits_{j \in S} {\left\langle {\psi _{i_{\pi  \left( j \right)} } } \right|A\left| {\psi _{i_j } } \right\rangle }  
= \left\langle {\psi _0 } \right|A\left| {\psi _0 } \right\rangle \left\langle {\psi _0 } \right|A\left| {\psi _0 } \right\rangle \left\langle {\psi _1 } \right|A\left| {\psi _1 } \right\rangle \geq 0
\]
and
\[
\prod\limits_{j \notin S} {\left\langle {\psi _{i_{\pi  \left( j \right)} } } \right|A\left| {\psi _{i_j } } \right\rangle }  
= \left\langle {\psi _0 } \right|A\left| {\psi _1 } \right\rangle \left\langle {\psi _1 } \right|A\left| {\psi _0 } \right\rangle \left\langle {\psi _0 } \right|A\left| {\psi _1 } \right\rangle \left\langle {\psi _1 } \right|A\left| {\psi _0 } \right\rangle  = \left| {\left\langle {\psi _0 } \right|A\left| {\psi _1 } \right\rangle } \right|^4 \geq 0. 
\]
Thus we have
\[
\prod\limits_{j = 1}^7 {\left\langle {\psi _{i_{\pi  \left( j \right)} } } \right|A\left| {\psi _{i_j } } \right\rangle }  
= \prod\limits_{j \in S} {\left\langle {\psi _{i_{\pi  \left( j \right)} } } \right|A\left| {\psi _{i_j } } \right\rangle }  \cdot \prod\limits_{j \notin S} {\left\langle {\psi _{i_{\pi  \left( j \right)} } } \right|A\left| {\psi _{i_j } } \right\rangle }  \ge 0.
\]
\end{Ex}

\begin{Rem}
For $T,A\in M_+(n,\mathbb{C})$, there exist $T,A$ and $p_1,p_2,\cdots,p_m$ 
such that 
$$Tr\left[ \prod _{i=1}^m \left(T^{p_i}A\right)  \right]  \notin \mathbb{R},$$
if $n \geq 3$ and $m\geq 3$.
For example, if we take \[
T = \left( \begin{array}{l}
 1\,\,\,\,0\,\,\,\,0 \\ 
 0\,\,\,\,2\,\,\,\,0 \\ 
 0\,\,\,\,0\,\,\,\,3 \\ 
 \end{array} \right),A = \left( \begin{array}{l}
 \,\,2\,\,\,\,\,\,\,\,\,\,\,i\,\,\,\,\,\,\,\,\,i \\ 
  -i\,\,\,\,\,\,\,\,\,2\,\,\,\,\,\,\,\,\,i \\ 
  -i\,\,\,\,\, -i\,\,\,\,\,\,2 \\ 
 \end{array} \right)
\]
and $p_1=1/6,p_2=1/3,p_3=1/2$, then $Tr\left[ \prod _{i=1}^3 \left(T^{p_i}A\right)  \right]$ approximately takes the complex number value $116.037+0.00260306i$.
Therefore the inequalities (\ref{ineq03}) does not make a sense in such more general cases than Theorem \ref{the01}.
\end{Rem}

We still have the following conjectures.
\begin{Con}
Do the following inequalities hold or not, for $T,A\in M_+(n,\mathbb{C})$ and positive numbers $p_1,p_2,\cdots,p_m$ with $p_1+p_2+\cdots +p_m=1$?
\begin{itemize}
\item[(i)] $Tr\left[ \left( T^{1/m}A\right)^m \right]\leq Re\left\{Tr\left[ T^{p_1}AT^{p_2}A\cdots T^{p_m}A \right]\right\}.$
\item[(ii)] $\vert Tr\left[ T^{p_1}AT^{p_2}A\cdots T^{p_m}A \right]\vert \leq Tr\left[ TA^m \right].$
\end{itemize}
\end{Con}

\section*{Acknowledgement}
The author (S.F.) was partially supported by the Japanese Ministry of Education, Science, Sports and Culture, 
Grant-in-Aid for Encouragement of Young Scientists (B), 20740067.
The authors (S.F. and K.Y.) was also partially supported by the Ministry of Education, Science, Sports and Culture, 
Grant-in-Aid for Scientific Research (B), 18300003.

\end{document}